\documentclass[12pt]{article}

\usepackage{amssymb}
\usepackage{epsfig}
\usepackage{rotating}

\def\b#1{\mbox{\boldmath $#1$}}
\def\bl#1{\mbox{\footnotesize \boldmath {$#1$}}} 

\newcommand{\logit}{{\rm logit}}
\newcommand{\expit}{{\rm expit}}

\newcommand{\se}{{\rm se}}
\newcommand{\pa}{\partial}         

\newcommand{\tr}{^{\prime}}

\newcommand{\al}{\alpha}
\newcommand{\be}{\beta}

\newcommand{\de}{\delta}

\newcommand{\Si}{\Sigma}
\newcommand{\la}{\lambda}
\renewcommand{\th}{\theta}

\newcommand{\convp}{\stackrel{p}\rightarrow}
\newcommand{\convd}{\stackrel{d}\rightarrow}

\newcommand{\ind}{\:\begin{sideways} $\hspace{-.4mm}\models$
\end{sideways}\:}

\begin{document}

\title{On the conditional logistic estimator\\ for repeated
binary outcomes\\ in two-arm experimental studies with
non-compliance}
\author{Francesco Bartolucci\footnote{ Department of Economics, Finance
and Statistics, University of Perugia, 06123 Perugia, Italy, {\em
email}: bart@stat.unipg.it}} \maketitle
\begin{abstract}
\noindent The behavior of the conditional logistic estimator is
analyzed under a causal model for two-arm experimental studies with
possible non-compliance in which the effect of the treatment is
measured by a binary response variable. We show that, when
non-compliance may only be observed in the treatment arm, the effect
(measured on the logit scale) of the treatment on compliers and that
of the control on non-compliers can be identified and consistently
estimated under mild conditions. The same does not happen for the
effect of the control on compliers. A simple correction of the
conditional logistic estimator is then proposed which allows us to
considerably reduce its bias in estimating this quantity and the
causal effect of the treatment over control on compliers. A two-step
estimator results whose asymptotic properties are studied by
exploiting the general theory on maximum likelihood estimation of
misspecified models. Finite-sample properties of the estimator are
studied by simulation and the extension to the case of missing
responses is outlined. The approach is illustrated by an application
to a dataset deriving from a study on the efficacy of a training
course on the practise of breast self examination.\vspace*{0.5cm}

\end{abstract}

\newpage
\section{Introduction}
Conditional logistic regression is a commonly used tool of data
analysis in the health sciences and medical statistics when the
outcome of interest is binary and subjects are observed before and
after a certain treatment or these subjects are somehow matched;
see, for instance, \cite{Breslow:1980}, \cite{Collett:1991},
\cite{Roth:1998} and \cite{Hosmer:2000}. The main reasons of the
popularity of the method are represented by its simplicity and by
the possibility of obtaining reliable estimates of the quantities of
interest under very mild conditions.

The first aim of this paper is that of illustrating the behavior of
the conditional logistic estimator when data come from two-arm
experimental studies with all-or-nothing compliance in which the
efficacy of the treatment is observed through a binary variable, and
a pre-treatment version of the same variable is available. We are
then in a context of repeated binary outcomes at two occasions,
before and after the treatment (or control), for which
non-compliance may represent a strong source of confounding in
estimating the causal effect of the treatment over control. An
example is represented by the study described by \cite{Ferro:1996}
on the effect of a training course on the attitude to practise
breast self examination (BSE); see also \cite{Mealli:2004}. In this
study, a significant number of women, among those randomly assigned
to the treatment, did not comply and preferred learning BSE by a
standard method (control). Moreover, the efficacy of the treatment
is measured by a binary variable indicating if a women regularly
practises BSE and another binary variable indicating if BSE is
practised in the proper way (provided it is practised).
Pre-treatment values of these response variables are also available.

In order to study the behavior of the conditional logistic estimator
in experimental studies of the type described above, we introduce a
causal model which includes parameters for the control and treatment
effects in the subpopulation of {\em compliers} and {\em
never-takers} and assumes that only the subjects assigned to the
treatment can access it. Given the type of experimental studies, we
assume that {\em always-takers} do not exist; we also assume that
{\em defiers} are no present. The model is then related to the
causal models described by \cite{Angrist:1996} and
\cite{Abadie:2003}; see also \cite{Rubin:2005}. It also allows for
the inclusion of base-line observable and unobservable covariates
which affect the response variables at the first and second
occasions. We show that, under this model, the conditional logistic
method allows us to identify and consistently estimate the effect of
the treatment on compliers and that of the control on never-takers.
However, apart from very particular cases, this method does not
allow us to identify the effect of control on the subpopulation of
compliers and then the causal effect of the treatment over control
on this subpopulation. As in other approaches for causal inference,
this effect is here measured on the logit scale; see
\cite{Ten:2003}, \cite{Van:2003}, \cite{Robins:2004} and
\cite{Laan:2007}.

Based on an approximation of the distribution of the observable
variables under the causal model, we then propose a correction for
the conditional logistic estimator which allows us to remove most of
its bias in estimating the effect of the treatment on compliers. It
results a two-step estimator which has some connection with the
estimator usually adopted for the selection model
\cite{Heckman:1979} and that proposed by \cite{Nag:2000} to estimate
the causal effect of a treatment in a context similar to the present
one. At the first step, the parameters of a model for the
probability that a subject is a complier are estimated. At the
second step, a conditional logistic likelihood is maximized which is
based on an approximated version of the conditional probability of
the response variables at the two occasions, given their sum. This
likelihood is computed on the basis of the first step parameter
estimates. The proposed estimator is very simple to use and is
consistent when the control has the same effect on compliers and
never-takers. This result holds regardless of the model that we
choose for the probability to comply. In the general case in which
compliers and never-takers react differently to control, the
estimator is not consistent but we show that it may converge in
probability to a value surprisingly close to the true value of the
causal parameters as the sample size grows to infinity. We also
derive a sandwich formula for its standard error. As we show, with
minor adjustments the two-step estimator leads to valid inference
even with missing responses.

The paper is organized as follows. In Section 2 we introduce the
causal model for repeated outcomes coming from two-arm experimental
studies. The behavior of the conditional logistic estimator under
this model is studied in Section 3. The correction of this estimator
is proposed in Section 4 where we also study the asymptotic and
finite-sample properties of the resulting two-step estimator. In
Section 5 we outline the extension of the approach to missing
responses and in Section 6 we provide an illustration based on an
application to the dataset deriving from the BSE study described
above. Final conclusions are reported in Section 7 where possible
extensions are also mentioned, such us that to experimental studies
in which subjects in both arms can access the treatment and then
non-compliance phenomena can be observed for all subjects.
\section{The causal model}\label{sec:ass}
Let $Y_1$ and $Y_2$ denote the binary response variables of
interest, let $\b V$ be a vector of observable covariates, let $Z$
be a binary variable equal to 1 when a subject is assigned to the
treatment and to 0 when he/she is assigned to the control and let
$X$ be the corresponding binary variable for the treatment actually
received. We recall that $\b V$ and $Y_1$ are pre-treatment
variables, whereas $Y_2$ is a post-treatment variable.
Non-compliance of the subjects involved in the experimental study
implies that $X$ may differ from $Z$. In particular we consider
experiments in which only subjects randomized to the treatment can
access it and therefore $Z=0$ implies $X=0$, whereas with $Z=1$ we
may observe either $X=0$ or $X=1$. Using a terminology taken from
\cite{Angrist:1996}, in this case we have only randomized
eligibility and we then consider two subpopulations: {\em compliers}
and {\em never-takers}. Nevertheless, the approach can be extended
to randomized experiments in which subjects in both arms can access
the treatment and therefore any configuration of $(Z,X)$ may be
observed; see Section \ref{sec:conc} for a discussion on this point.
In both types of experiment, we assume that defiers are not present
and our aim is that of estimating the causal effect of the treatment
over control in the subpopulation of compliers.

We assume that the behaviour of a subject depends on the observable
covariates $\b V$, a latent variable $U$ representing the effect of
unobservable covariates on both response variables and a latent
variable $C$ representing the attitude to comply with the assigned
treatment. The last one, in particular, is a discrete variable with
two levels: 0 for never-takers, 1 for compliers.

The model is based on the following assumptions:
\begin{itemize}
\item[A1:] $C\ind Y_1|(U,\b V)$, i.e. $C$ is conditionally
independent of $Y_1$ given $(U,\b V)$;
\item[A2:] $Z\ind(U,Y_1,C)|\b V$;
\item[A3:] $X\ind(U,\b V,Y_1)|(C,Z)$ and, with probability 1, $X=Z$
when $C=1$ (compliers) and $X=0$ when $C=0$ (never-takers);
\item[A4:] $Y_2\ind (Y_1,Z)|(U,\b V,C,X)$;
\item[A5:] for any $u,\b v,c$ and $x$, we have
\begin{equation}
\logit[p(Y_2=1|u,\b v,c,x)]-\logit[p(Y_1=1|u,\b v)]=\b a(\b
v,x)\tr\b\al+c(1-x)\b b(\b v)\tr\b\be,\label{ass5}
\end{equation}
where $\b a(\b v,x)$ and $\b b(\b v)$ are known functions which
depend on observable covariates and (limited to the first) on the
received treatment and $\b\al$ and $\b\be$ are corresponding vectors
of parameters.
\end{itemize}

The above assumptions lead to a dependence structure on the
observable and unobservable variables which is represented by the
DAG in Figure 1.

\begin{figure}[ht]\centering
\epsfig{figure=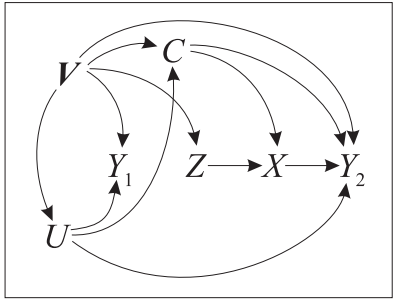,height=4cm} \caption{\em DAG for the
model based on assumptions A1-A5. $U$ and $\b V$ represent
unobservable and observable covariates affecting the response
variables $Y_1$ and $Y_2$, $C$ is a binary variable for the
compliance status and $Z$ and $X$ are binary variables for the
assigned and received treatment.}\label{fig1}\vspace*{0.25cm}
\end{figure}

Assumption A1 says that the tendency to comply only depends on
$(U,\b V)$. Assumption A2 is satisfied in randomized experiments,
even when randomization is conditioned on the observable covariates.
This assumption could be relaxed by requiring that $Z$ is
conditionally independent of $U$ given $(\b V,Y_1)$, so that
randomization can also be conditioned on the first outcome.
Assumption A3 is rather obvious considering that $C$ represents the
tendency of a subject to comply with the assigned treatment.
Assumption A4 implies that there is no direct effect of $Y_1$ on
$Y_2$, since the distribution of the latter only depends on $(U,\b
V,C,X)$. Using a terminology which is well known in the literature
on latent variable models, this is an assumption of {\em local
independence}. Assumption A4 also implies an assumption known in the
causal inference literature as {\em exclusion restriction},
according to which $Z$ affects $Y_2$ only through $X$. Finally,
assumption A5 says that the distribution of $Y_2$ depends on the
vectors of parameters $\b\al$ and $\b\be$ through the functions $\b
a(\b v,x)$ and $\b b(\b v)$. Note, in particular, that $\b a(\b
v,0)\tr\b\al$ is the effect of the control on never-takers, $\b a(\b
v,1)\tr\b\al$ is the effect of the treatment on compliers, whereas
$\b b(\b v)\tr\b\be$ is the differential effect of the control
between compliers and never-takers. In the simplest case, we have
\begin{equation}
\b a(\b v,x) = (1-x,x)\tr\quad\mbox{and}\quad b(\b v) = 1,
\label{par2}
\end{equation}
so that $\b\al = (\al_1,\al_2)\tr$ and $\be$ have an obvious
interpretation as effects of the control and the treatment on
specific subpopulations.

As mentioned above, the most interest quantity to estimate is the
{\em causal effect} of the treatment over the control in the
subpopulation of compliers. In the present approach, this effect is
defined as the difference in logits (log-odds ratio)
\begin{eqnarray}
\de(\b v)&=&\logit[p(Y_2=1|u,\b v,C=1,X=1)]-\logit[p(Y_2=1|u,\b
v,C=1,X=0)]=\nonumber\\
&=&[\b a(\b v,1)-\b a(\b v,0)]\tr\b\al-\b b(\b
v)\tr\b\be,\label{eq:caus_eff}
\end{eqnarray}
and, with reference to the subpopulation of compliers, it
corresponds to the increase of the logit of the probability of
success when $X$ goes from 0 to 1, all the other factors remaining
unchanged. This quantity depends on the covariates in $\b V$ and
then an overall causal effect can be computed as the average of
$\de(\b v)$ over suitable configurations of these covariates. Under
(\ref{par2}), we simply have $\de=\al_2-\al_1-\be$ and then
computing this average is not necessary. Also note that the model
makes sense, not only when $Y_1$ is a response variable of the same
nature of $Y_2$ that is observed before the treatment, but also when
$Y_1$ is a variable which is affected neither by the compliance
status nor by the treatment received and such that the difference
between the logits in (\ref{ass5}) is independent of $u$ and $\b v$.

That based on assumptions A1-A5 is a causal model in the sense of
\cite{Pearl:1995} since all the observable and unobservable factors
affecting the response variables of interest are included. Indeed,
the same model could be formulated by exploiting potential outcomes,
which we denote by $Y_2^{(z,x)}$, $z,x=0,1$. In this case, the model
could be formulated on the basis of assumptions A1-A3 and the
following assumptions which substitute A4-A5:
\begin{itemize}
\item[A4$^*$:] $Y_2^{(z,x)}=Y_2^{(x)}$ for $z,x=0,1$ ({\em exclusion
restriction}) and $(Y_2^{(0)},Y_2^{(1)})\ind(Y_1,Z,X)|(U,\b V,C)$;
\item[A5$^*$:] for any $u,\b v,c$ and $x$, we have
\begin{equation}
\logit[p(Y_2^{(x)}=1|u,\b v,c)]-\logit[p(Y_1=1|u,\b v)]=\b a(\b
v,x)\tr\b\al+c(1-x)\b b(\b v)\tr\b\be;\label{potential}
\end{equation}
\item[A6$^*$:] $Y_2=Y_2^{(x)}$ for any given value of $U,\b V,C$ and
any given value $x$ of $X$.
\end{itemize}
It may be easily realized that the model based on assumptions A1-A5
is equivalent to that based on assumptions A1-A3 and A4$^*$-A6$^*$.
In a similar context, this kind of equivalence between causal models
is dealt with by \cite{Ten:2003}.

It is even more clear from (\ref{potential}) that $\b\al$ and
$\b\be$ are vectors of causal parameters, in the sense that they
allow us to measure the causal effect of the treatment over control
in the subpopulation of compliers. Using potential outcomes and for
simplicity under (\ref{par2}), this effect may now be expressed as
\[
\de=\logit[p(Y_2^{(1)}=1|u,\b v,C=1)]-\logit[p(Y_2^{(0)}=1|u,\b
v,C=1)]=\al_2-\al_1-\be.
\]

\section{Behavior of the conditional logistic
estimator}\label{sec:cond}
In this section we show that, by the conditional logistic method and
under mild conditions, we can identify and consistently estimate the
parameter vector $\b\al$ which measures the effect of the treatment
on compliers and that of the control on never-takers. However, the
same is not possible for $\b\be$. In order to show this we need to
derive the conditional distribution of $(Y_1,Z,X,Y_2)$ given $\b V$
under the assumptions introduced in Section \ref{sec:ass}. The
probability mass function of this distribution is denoted by
$f(y_1,z,x,y_2|\b v)$.

First of all, these assumptions imply that the probability function
of the conditional distribution of $(Y_1,Z,X,Y_2)$ given $(U,\b V,C)$
is equal to
\[
p(y_1,z,x,y_2|u,\b v,c)=p(y_1|u,\b v)p(z|\b v)p(x|c,z)p(y_2|u,\b
v,c,x).
\]
After some algebra we have
\[
p(y_1,z,x,y_2|u,\b v)=\frac{e^{(y_1+y_2)\la(u,\bl
v)}}{1+e^{\la(u,\bl v)}}p(z|\bl v)\sum_c p(x|c,z)\frac{e^{y_2t(c,\bl
v,x)}}{1+e^{\la(u,\bl v)+t(c,\bl v,x)}}\pi(c|u,\b v),
\]
where $\la(u,\b v)=\logit[p(Y_1=1|u,\b v)]$, $\pi(c|u,\b v)=p(c|u,\b
v)$ and the sum $\sum_c$ is extended to $c=0,1$. Moreover, $t(c,\b
v,x)=\b a(\b v,x)\tr\b\al+c(1-x)\b b(\b v)\tr\b\be$; see equation
(\ref{ass5}). Finally, denoting by $\phi(u|\b v)$ the density
function of the conditional distribution of $U$ given $\b V$ and
letting $y_+=y_1+y_2$, we have
\begin{equation}
f(y_1,z,x,y_2|\b v)=p(z|\b v)\int\frac{e^{y_+\la(u,\bl
v)}}{1+e^{\la(u,\bl v)}}\sum_c p(x|c,z)\frac{e^{y_2t(c,\bl
v,x)}}{1+e^{\la(u,\bl v)+t(c,\bl v,x)}}\pi(c|u,\b v)\phi(u|\b
v)du.\label{p_man}
\end{equation}

When $z=1$ (treatment arm), $p(x|c,z)$ is equal to 1 for $x=c$ and to
0 otherwise. Consequently, (\ref{p_man}) reduces to
\begin{equation}
f(y_1,1,x,y_2|\b v)=p(Z=1|\b v)\int\frac{e^{y_+\la(u,\bl
v)}}{1+e^{\la(u,\bl v)}} \frac{e^{y_2\bl a(\bl
v,x)\tr\bl\al}}{1+e^{\la(u,\bl v)+\bl a(\bl
v,x)\tr\bl\al}}\pi(x|u,\bl v)\phi(u|v)du,\label{p_man1}
\end{equation}
since it is possible to identify a subject in this arm as a
never-taker or a compiler according to whether $x=0$ or $x=1$ and
$t(c,\b v,x)=\b a(\b v,x)\tr\b\al$ when $c=x$. Now let $f(y_2|\b
v,z,x)$ denote the probability mass function of the conditional
distribution of $Y_2$ given $(\b V,X,Y_+=1)$, with $Y_+=Y_1+Y_2$.
The above result implies that
\begin{equation}
f(y_2|\b v,1,x)=\frac{f(1-y_2,1,x,y_2|v)}{f(0,1,x,1|v)+
f(1,1,x,0|v)}=\frac{e^{y_2\bl a(\bl v,x)\tr\bl\al}}{1+e^{\bl a(\bl
v,x)\tr\bl\al}}.\label{eq:pcond}
\end{equation}
This probability function depends neither on the distribution of $U$
nor on the function $\la(u,\b v)$. Consequently, as already
mentioned above, we can identify the parameters in $\b\al$ measuring
the effect of the treatment on compliers and that of the control on
never-takers and the conditional logistic estimator of these
parameters is consistent. It is also worth to observe that under
assumption (\ref{par2}), which implies that
$\b\al=(\al_1,\al_2)\tr$, the conditional logistic estimator of
these parameters has an explicit form given by
\[
\hat{\al}_1 = \log\frac{n_{0101}}{n_{1100}}\quad\mbox{and}\quad
\hat{\al}_2 = \log\frac{n_{0111}}{n_{1110}},
\]
where $n_{y_1zxy_2}$ denotes the number of subjects in the sample
with response configuration $(y_1,y_2)$ who are in the control or
treatment arm (according to whether $z=0$ or $z=1$) and chose or not
the treatment (according to whether $x=0$ or $x=1$).

When $z=0$ (control arm), $p(x|c,z)$ is equal to 1 for $x=0$ and to 0
otherwise (regardless of $c$), since no subject in this arm can
access the treatment. Then, expression (\ref{p_man}) reduces to
\begin{eqnarray*}
f(y_1,0,0,y_2|\b v)&=&p(Z=0|\b v)g(y_1,y_2|\b v),\nonumber\\
g(y_1,y_2|\b v)&=&\int\frac{e^{y_+\la(u,\bl v)}}{1+e^{\la(u,\bl
v)}}\sum_c \frac{e^{y_2t(c,\bl v,0)}}{1+e^{\la(u,\bl v)+t(c,\bl
v,0)}}\pi(c|u,\b v)\phi(u|\b v)du,
\end{eqnarray*}
which is more complex than (\ref{p_man1}), being based on a mixture
between the conditional distribution of $Y_2$ for the subpopulation
of compliers and for that of never-takers. In this case, we cannot
remove the dependence on the distribution of $U$ and on $\la(u,\b
v)$ via conditioning on $Y_+$. Then, we cannot identify and
consistently estimate the parameters in $\b\be$ corresponding to the
differential effect of control on compliers with respect to
never-takers. The same happens for the causal effect $\de(\b v)$
defined in (\ref{eq:caus_eff}).

In order to better investigate this point we consider an
approximation of $\log g(y_1,y_2|\b v)$ based on a first-order
Taylor series expansion around $\b\be=\b 0$, with $\b 0$ denoting a
column vector of zeros of suitable dimension. Note that this point
corresponds to the situation in which the control has the same
effect on compliers and never-takers. We have that
\[
\log g(y_1,y_2|\b v)\approx\log g_0(y_1,y_2|\b v)+h(y_1,y_2|\b v)\b
b(\b v)\tr\b\be,
\]
where $g_0(y_1,y_2|\b v)$ is equal to $g(y_1,y_2|\b v)$ computed at
$\b\be=\b 0$, that is
\[
g_0(y_1,y_2|\b v)=\int\frac{e^{y_+\la(u,\bl v)}}{1+e^{\la(u,\bl
v)}}\frac{e^{y_2t(0,\bl v,0)}}{1+e^{\la(u,\bl v)+t(0,\bl
v,0)}}\phi(u|\b v)du,
\]
and
\[
h(y_1,y_2|\b v)=\frac{2y_2-1}{g_0(y_1,y_2|\b
v)}\int\frac{e^{y_+\la(u,\bl v)}}{1+e^{\la(u,\bl
v)}}\frac{e^{t(0,\bl v,0)}}{[1+e^{\la(u,\bl v)+t(0,\bl
v,0)}]^2}\pi(1|u,\b v)\phi(u|\b v)du.
\]
Now, because $g_0(0,1|\b v)=g_0(1,0|\b v)e^{t(0,\bl v,0)}$ and
recalling that $t(0,\b v,0) = \b a(\b v,0)\tr\b\al$, after some
algebra (see Appendix A1 for details) we find
\begin{equation}
\log\frac{f(0,0,0,1|\b v)}{f(1,0,0,0|\b v)}=\log g(0,1|\b v)-\log
g(1,0|\b v)\approx\b a(\b v,0)\tr\b\al+h(\b v)\b b(\b v)\tr\b\be,
\label{eq:approx1}
\end{equation}
where $h(\b v)$ is a correction factor defined as
\begin{equation}
h(\b v)=\frac{1}{g_0(1,0|\b v)}\int\frac{e^{\la(u,\bl
v)}}{1+e^{\la(u,\bl v)}}\frac{1}{1+e^{\la(u,\bl v)+t(0,\bl
v,0)}}\pi(1|u,\b v)\phi(u|\b v)du.\label{eq:r}
\end{equation}
This correction term is simply equal to $\pi(1|\b v)=p(C=1|\b v)$
when $C$ is conditionally independent of $U$ given $\b V$. We then
have
\begin{equation}
f(y_2|\b v,0,0)=\frac{f(1-y_2,0,0,y_2|\b v)}{f(0,0,0,1|\b
v)+f(1,0,0,0|\b v)}\approx\frac{e^{\bl a(\bl v,0)\tr\bl\al+h(\bl
v)\bl b(\bl v)\tr\bl\be}}{1+e^{\bl a(\bl v,0)\tr\bl\al+h(\bl v)\bl
b(\bl v)\tr\bl\be}}\label{eq:approx}
\end{equation}
which shows that a conditional logistic estimator based on
regressing $Y_2$ on $\b V$, only for the cases in which $Y_+=1$ and
$Z=0$, would estimate a quantity which does not correspond to the
effect of the control either on compliers or on never-takers.

To clarify the above point consider the case of absence of
covariates in which assumption (\ref{par2}) holds. In this case, the
conditional logistic estimator is equal to $\log(n_{0001}/n_{1000})$
which converges in probability to a quantity close to $\al_1+h\be$.
Then, provided that $h$ can be suitably estimated, a correction for
this estimator can be implemented so as to reduce its bias. This
idea is exploited to propose a general approach which allows us to
considerably reduce the bias of the conditional logistic estimator
applied to the data coming from the control group.
\section{Corrected conditional logistic estimator}\label{sec:corr}
With reference to a sample of $n$ subjects included in the two-arm
experimental study, let $y_{i1}$ denote the observed value of $Y_1$
for subject $i$, let $y_{i2}$ denote the value of $Y_2$ for the same
subject and let $\b v_i$, $z_i$ and $x_i$ denote the corresponding
values of $\b V$, $Z$ and $X$, respectively.
\subsection{The estimator}\label{sec:corr_est}
In order to estimate the parameters of the causal model introduced
in Section \ref{sec:ass}, we rely on a standard logistic regression
applied to the data coming from the treatment arm (for which we can
disentangle compliers from never-takers) and a logistic regression
based on approximation (\ref{eq:approx}) for the data coming from
the control arm. Note that to exploit this approximation we need to
estimate the correction term $h(\b v)$ defined in (\ref{eq:r}). For
sake of simplicity, we assume that $C$ is conditionally independent
on $U$ given $\b V$, so that this correction term corresponds to
$\pi(1|\b v)$ and for the latter we assume the logit model
\begin{equation}
\logit[\pi(1|\b v)] = \b m(\b v)\tr\b\eta,\label{logit}
\end{equation}
where $\b m(\b v)$ is a known function of the observed covariates.
The implications of this assumption will be studied in the following
(see Section \ref{sec:corr_est_asy}). It results an estimator of the
causal effect parameters whose main advantage is the simplicity of
use. The estimator recalls the two-step estimator of the selection
model \cite{Heckman:1979} and the estimator proposed by
\cite{Nag:2000} in the causal inference literature.

The two steps on which the proposed estimation method is based are
the following:
\begin{enumerate}
\item {\em Estimation of $\b\eta$}. Since the
compliance status may be directly observed for those subjects
assigned to the treatment arm, estimation of $\b\eta$ is based on
the observed values $x_i$ and $\b v_i$ for every $i$ such that
$z_i=1$. Taking into account that the distribution of $Z$ is allowed
to depend on $\b V$, we then proceed by maximizing the weighted
log-likelihood
\[
\ell_1(\b\eta)=\sum_i\frac{z_i}{p(z_i|\b v_i)}[x_i\log\pi(1|\b
v_i)+(1-x_i)\log\pi(0|\b v_i)],
\]
with weights corresponding to the inverse probabilities $1/p(z_i|\b
v_i)$.
\item {\em Estimation of $\b\al$ and $\b\be$}. This is done by maximizing the conditional
log-likelihood
\begin{equation}
\ell_2(\b\al,\b\be;\b\eta) = \sum_i
d_i\ell_{i2}(\b\al,\b\be;\b\eta),\label{lk2}
\end{equation}
where $\ell_{i2}(\b\al,\b\be;\b\eta)$ corresponds to the logarithm
of (\ref{eq:pcond}) when $z_i=1$ and to that of (\ref{eq:approx})
when $z_i=0$, once the parameter vector $\b\eta$ has been
substituted by its estimate $\hat{\b\eta}$ obtained at the first
step. Moreover, $d_i$ is a dummy variable equal to 1 if $y_{i+}=1$,
with $y_{i+}=y_{i1}+y_{i2}$, and to 0 otherwise, so that subjects
with response configuration $(0,0)$ or $(1,1)$ are excluded since
the conditional probability of these configurations given their sum
would not depend either on $\b\al$ or $\b\be$.
\end{enumerate}

Maximization of $\ell_1(\b\eta)$ may be performed by a standard
Newton-type algorithm; the data to be used are only those concerning
the subjects in the treatment arm. The same algorithm may used to
maximize $\ell_2(\b\al,\b\be;\b\eta)$. In this case, the data to be
used concern all subjects with sum of the response variables (before
and after) equal to 1. Moreover, collecting the parameter vectors
$\b\al$ and $\b\be$ in a unique vector
$\b\th=(\b\al\tr,\b\be\tr)\tr$, the design matrix to be used has
rows $\b w(\b v_i,z_i,x_i)\tr$, where
\[
\b w(\b v_i,z_i,x_i)=\pmatrix{\b a(\b v_i,0)\cr (1-z_i)\pi(1|\b
v_i)\b b(\b x_i)},
\]
which corresponds to $(\b a(\b v_i,x_i)\tr,\b 0\tr)$ when $z_i=1$
and $(\b a(\b v_i,0)\tr,\pi(1|\b v_i)\b b(\b x_i)\tr)$ when $z_i=0$.
At the end, by substituting the subvectors $\hat{\b\al}$ and
$\hat{\b\be}$ of $\hat{\b\th}$ into (\ref{eq:caus_eff}) we obtain an
estimate $\hat{\de}(\b v)$ of the causal effect of the treatment,
with respect to control, for compliers with covariate configuration
$\b v$. Under (\ref{par2}), this estimate reduces to
$\hat{\de}=\hat{\al}_2-\hat{\al}_1-\hat{\be}$.

With small samples, which are not uncommon in certain experimental
studies, it might happen that discordant response configurations, of
type (0,1) or (1,0), are not observed for certain configurations of
$(Z,X)$. This would imply that the estimator of $\b\th$ cannot be
computed. To overcome this problem, we follow a {\em rule of thumb}
consisting of: ({\em i}) checking that both discordant response
configurations are present for each observable configuration of
$(X,Z)$; ({\em ii}) adding to the dataset the response
configurations which are missing. To the added response
configurations we assign a vector of covariates $\b v_i$ equal to
the sample average and weight 0.5 in the conditional log-likelihood
(\ref{lk2}). The simulation study in Section 4.2.2 allows us to
evaluate the impact of this correction on the inferential properties
of the estimator.

Concerning the estimation of the variance-covariance matrix of the
estimators $\hat{\b\al}$ and $\hat{\b\be}$, consider that
$(\hat{\b\eta}\tr,\hat{\b\th}\tr)\tr$ correspond to the solution
with respect to $(\b\eta\tr,\b\th)\tr$ of the equation $\b
s(\b\eta,\b\th)=\b 0$, where
\[
\b
s(\b\eta,\b\th)=\left(\begin{array}{c}{\displaystyle\frac{\pa\ell_1(\b\eta)}{\pa\b\eta}}\vspace*{0.1cm}\cr
{\displaystyle\frac{\pa\ell_2(\b\th;\b\eta)}{\pa\b\th}}\end{array}\right),
\]
with $\ell_2(\b\th;\b\eta)=\ell_2(\b\al,\b\be;\b\eta)$. When logit
model (\ref{logit}) is assumed, the first subvector of $\b
s(\b\eta,\b\th)$ is equal to
\[
\frac{\pa\ell_1(\b\eta)}{\pa\b\eta}=\sum_i\frac{z_i}{p(z_i|\b
v_i)}[x_i-\pi(1|\b v_i)]\b m(\b v_i),
\]
whereas the second subvector is equal to
\[
\frac{\pa\ell_2(\b\th)}{\pa\b\th}=\sum_id_i\bigg[y_{i2}-\frac{e^{\bl
w(\bl v_i,z_i,x_i)\tr\bl\th}}{1+e^{\bl w(\bl
v_i,z_i,x_i)\tr\bl\th}}\bigg]\b w(\b v_i,z_i,x_i).
\]
From \cite{Huber:1967} and \cite{White:1982}, the following sandwich
estimator of the variance-covariance matrix of
$(\hat{\b\eta}\tr,\hat{\b\th}\tr)\tr$ results
\begin{equation}
\hat{\b\Si}(\hat{\b\eta},\hat{\b\th}) = \hat{\b H}^{-1}\hat{\b
K}(\hat{\b H\tr})^{-1},\label{sandwich}
\end{equation}
where $\hat{\b H}$ is the derivative of $\b s(\b\eta,\b\th)$ with
respect to $(\b\eta\tr,\b\th\tr)$ and $\hat{\b K}$ is an estimate of
the variance-covariance matrix of $\b s(\b\eta,\b\th)$, both
computed at $\b\eta=\hat{\b\eta}$ and $\b\th=\hat{\b\th}$. Explicit
expressions for these matrices are given in Appendix A2. We can then
obtain an estimate of the variance-covariance matrix of
$\hat{\b\th}$, denoted by $\hat{\b\Si}(\hat{\b\th})$ or
alternatively by $\hat{\b\Si}(\hat{\b\al},\hat{\b\be})$, as a
suitable block of the matrix
$\hat{\b\Si}(\hat{\b\eta},\hat{\b\th})$. We can also obtain the
standard error for estimator of the causal effect $\hat{\de}(\b v)$.
In particular, when (\ref{par2}) holds then the standard error for
$\hat{\de}$ is simply
$\se(\hat{\de})=\sqrt{\b\Delta\tr\hat{\b\Si}(\hat{\b\th})\b\Delta}$,
where $\b\Delta=(-1,1,-1)\tr$ is a vector such that
$\hat{\de}=\b\Delta\tr\hat{\b\th}$.
\subsection{Properties of the two-step estimator}
\subsubsection{Asymptotic properties}\label{sec:corr_est_asy}
Suppose that $\b v_i$, $y_{i1}$, $z_i$, $x_i$ and $y_{i2}$, with
$i=1,\ldots,n$, are independently drawn from the {\em true model}
based on assumptions A1-A5 with parameters $\b\al=\b\al_0$ and
$\b\be=\b\be_0$. This model must ensure that
\[
f(y_1,z,x,y_2|\b v)>0\quad\mbox{for all }\quad \b v,y_1,z,x,y_2.
\]
Provided that the functions $\b a(\b v,x)$ and $\b b(\b v)$ satisfy
standard regularity conditions, which are necessary to ensure that
the expected value of the second derivative of
$\ell_2(\b\al,\b\be;\b\eta)/n$ is of full rank, the theory on
maximum likelihood estimation of misspecified models of
\cite{White:1982}, see also \cite{Newey:1994}, implies that the
two-step estimators $\hat{\b\al}$ and $\hat{\b\be}$ satisfy the
following asymptotic properties as $n\rightarrow\infty$:
\begin{itemize}
\item {\em consistency}: $\hat{\b\al}\convp\b\al_*$ and
$\hat{\b\be}\convp\b\be_*$, with $\b\al_*$ and $\b\be_*$ being the
{\em pseudo-true parameter vectors} which are equal, respectively,
to true parameter vectors $\b\al_0$ and $\b\be_0$ when $\b\be_0=\b
0$;
\item {\em asymptotic normality}:
$\sqrt{n}(\hat{\b\al}\tr,\hat{\b\be}\tr)\tr\convd N[\b
0,\b\Omega(\b\al_*,\b\be_*)]$, with $\b\Omega(\b\al_*,\b\be_*)$
being the limit in probability of the matrix
$\hat{\b\Si}(\hat{\b\al},\hat{\b\be})/n$; see definition
(\ref{sandwich}).
\end{itemize}
We recall that $\convp$ means convergence in probability, whereas
$\convd$ means converge in distribution. Moreover, in order to give
a formal definition of $\b\al_*$ and $\b\be_*$ we have to consider
that these correspond to the supremum of
$E_0[\ell_2(\b\al,\b\be;\b\eta_*)/n]$, where $E_0$ denotes
expectation under the true model and $\b\eta_*$ denotes the limit in
probability of the estimator $\hat{\b\eta}$ computed at the first
step. Clearly, since the log-likelihood
$\ell_2(\b\al,\b\be;\b\eta_*)$ is based on an approximation of the
true model around $\b\be=\b 0$, we have that $\b\al_*=\b\al_0$ and
$\b\be_*=\b\be_0$ when $\b\be_0=\b 0$. This implies that in this
case $\hat{\de}(\b v)\convp\de_0(\b v)$, with $\de_0(\b v)$ being
equal to the true causal effect of the treatment over control for a
complier with covariates $\b v$. Obviously, this is not ensured when
$\b\be_0=\b 0$, but we expect $\b\al_*$ and $\b\be_*$ to be
reasonably close to, respectively, $\b\al_0$ and $\b\be_0$ when
$\b\be_0$ is not too far from $\b 0$. The same may be said about the
estimator $\hat{\de}(\b v)$ of $\de(\b v)$.

In order to illustrate the previous point, we considered a true
model which involves only one observable covariate $V$ and under
which the joint distribution of $(U,V)$ is
\begin{equation}
N\left[\pmatrix{0\cr 0},\pmatrix{1 & \rho\cr \rho &
1}\right],\label{simu1}
\end{equation}
with $\rho=0.00,0.75$. Moreover, we assumed (\ref{par2}) and that
$Y_1$, $C$, $Z$, $Y_2$ have Bernoulli distribution with
probabilities of success chosen, respectively, as follows
\begin{equation}
\begin{array}{rll}
p(Y_1=1|u,v)&=&\expit[(u+v)/\sqrt{1+\rho^2}-1],\vspace*{0.25cm}\cr
\pi(c|u,v)&=&\expit[(u+v)/\sqrt{1+\rho^2}/2],\vspace*{0.25cm}\cr
p(Z=1|v) &=&\expit(-v),\vspace*{0.25cm}\cr
p(Y_2=1|u,v,c,x)&=&\expit[(u+v)/\sqrt{1+\rho^2}-1+(1-x)\al_1+x\al_2+c(1-x)\be],\label{simu2}
\end{array}
\end{equation}
where $\expit(\cdot)$ is the inverse function of $\logit(\cdot)$,
$\be=-1.00,-0.75,\ldots,1.00$ and $\al_1$ is defined so that the
casual effect $\de=\al_2-\al_1-\be$ is equal to 0 or 1, with
$\al_2=1$ when $\de=0$ and $\al_2=2$ when $\de=1$. Under this model,
we computed the limit in probability $\de_*$ of each of the
following estimators:
\begin{itemize}
\item $\hat{\de}_{null}$: two-step conditional logistic estimator of
$\de$ in which the probability to comply is assumed to do not depend
on the covariate; this is equivalent to letting $m(v)=1$ in
(\ref{logit});
\item $\hat{\de}_{cov}$: as above with $\b m(v) = (1,v)\tr$, so that the
covariate is also used to predict the probability to comply;
\item $\hat{\de}_{itt}$: {\em intention to treat} (ITT) estimator
based on the conditional logistic regression of $Y_2$ on $(1,Z)$
given $Y_+=1$;
\item $\hat{\de}_{tr}$: {\em treatment received} (TR)
estimator based on the conditional logistic regression of $Y_2$ on
$(1,X)$  given $Y_+=1$.
\end{itemize}
The limit in probability of these estimators is represented, with
respect to the true value of $\be$, in Figure 2.

\begin{figure}[ht]\centering
\epsfig{figure=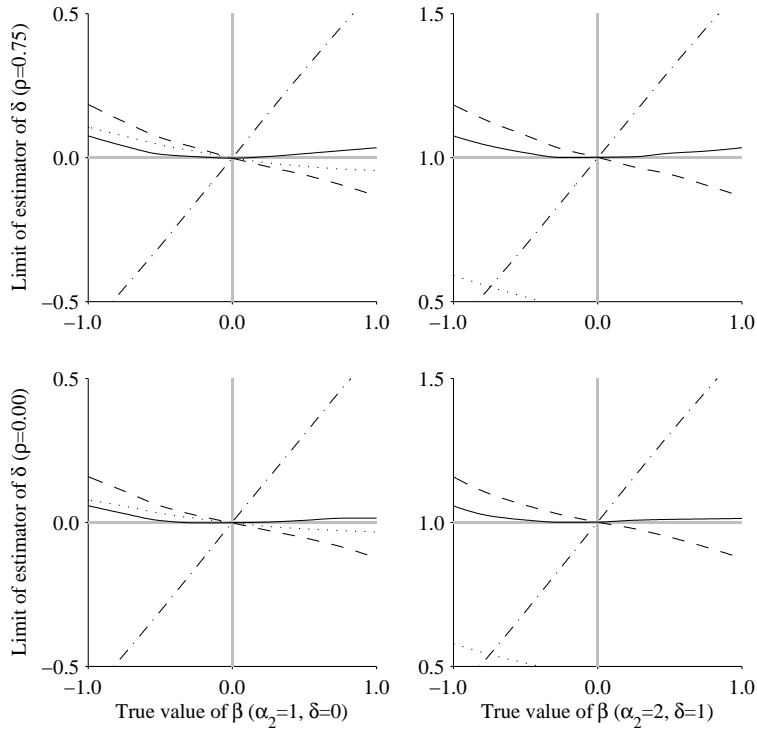,width=10cm} \caption{\em Limit of
the estimators $\hat{\de}_{null}$ (dashed line), $\hat{\de}_{cov}$
(solid line), $\hat{\de}_{itt}$ (dotted line) and $\hat{\de}_{tr}$
(dash-dotted line) under assumptions (\ref{simu1}) and
(\ref{simu2}), with $\rho=0.00,0.75$, $\be$ between -1 and 1 and
$\de=0,1$, with $\al_1$ and $\al_2$ defined
accordingly.}\label{fig2}\vspace*{0.25cm}
\end{figure}

It may be observed that, when the true value of $\be$ is 0, the
limit $\de_*$ is equal to true value of $\de$ for both estimators
$\hat{\de}_{null}$ and $\hat{\de}_{cov}$. This is in agreement with
our conclusion above about the asymptotic behavior of the proposed
estimator. When the true value of $\be$ is different from 0,
instead, this does not happen but, at least for $\hat{\de}_{cov}$,
the distance of $\de_*$ from the true $\de$ is surprisingly small
and does not seem to be affected by the correlation between $U$ and
$V$ which is measured by $\rho$. We recall that, although our
estimator is derived under the assumption that $C$ is conditionally
independent of $U$ given $\b V$, this result is obtained under a
model which assumes that both $U$ and $\b V$ have a direct effect on
$C$.

A final points concerns the ITT and TR estimators. The first is
adequate only if the true value of $\de$ is equal to 0 (plots on the
left of Figure 1), whereas it is completely inadequate when it is
equal to 1 (plots on the right). The TR estimator, instead, is
consistent only when the true value of $\be$ is equal to 0, but in
the other cases it has a strong bias. Overall, even if based on a
logistic regression method, these two estimators behave much worse
than our estimator.
\subsubsection{Finite-sample properties}
In order to assess the finite-sample properties of the two-step
estimator, we performed a simulation study based on 1000 samples of
size $n=200,500$ generated from the model based on assumptions
(\ref{simu1}) and (\ref{simu2}). For each simulated sample, we
computed the estimators $\hat{\b\al}=(\hat{\al}_1,\hat{\al}_2)\tr$,
$\hat{\be}$ and $\hat{\de}$ based of a model for the probability to
comply of type (\ref{logit}), with $\b m(v) = (1,v)\tr$. Using the
notation of Section \ref{sec:corr_est_asy}, these estimators could
also be denoted by $\hat{\b\al}_{cov}$, $\hat{\be}_{cov}$ and
$\hat{\de}_{cov}$. The results, in term of bias and standard
deviation of the estimators and in terms of mean of the standard
errors, are reported in Table \ref{tab1} (when the true value of
$\de$ is 0) and in Table \ref{tab2} (when the true value of $\de$ is
1). Note that, for small sample sizes as those we are considering
here, it may happen that there are not discordant configurations.
Consequently, we apply the rule of thumb described in Section
\ref{sec:corr_est} to prevent instabilities of the estimator.

\begin{table}[ht!]\centering
{\small\begin{tabular}{ccccclcccc}
\hline\hline        \vspace*{-0.4cm}                    \\
$n$ & $\rho$  &   $\al_1$   &   $\al_2$ &   $\be$ & & $\hat{\al}_1$
&   $\hat{\al}_2$   &   $\hat{\be}$   & $\hat{\de}$ \\
\hline
200 &   0.00    &   2   &   1   &   -1  &   bias    &   0.072   &   0.045   &   -0.144  &   0.117   \\
    &       &       &       &       &   st.dev. &   (0.606) &   (0.596) &   (1.344) &   (1.099) \\
    &       &       &       &       &   mean s.e.   &   (0.658) &   (0.573) &   (1.410) &   (1.097) \\\cline{3-10}
    &       &   1   &   1   &   0   &   bias    &   0.078   &   0.070   &   -0.024  &   0.017   \\
    &       &       &       &       &   st.dev. &   (0.575) &   (0.622) &   (1.268) &   (1.041) \\
    &       &       &       &       &   mean s.e.   &   (0.555) &   (0.578) &   (1.222) &   (1.019) \\\cline{3-10}
    &       &   0   &   1   &   1   &   bias    &   0.024   &   0.096   &   0.021   &   0.051   \\
    &       &       &       &       &   st.dev. &   (0.556) &   (0.614) &   (1.226) &   (1.037) \\
    &       &       &       &       &   mean s.e.   &   (0.553) &   (0.585) &   (1.204) &   (1.000) \\\hline
200 &   0.75    &   2   &   1   &   -1  &   bias    &   0.122   &   0.098   &   -0.178  &   0.153   \\
    &       &       &       &       &   st.dev. &   (0.621) &   (0.632) &   (1.386) &   (1.171) \\
    &       &       &       &       &   mean s.e.   &   (0.658) &   (0.602) &   (1.380) &   (1.093) \\\cline{3-10}
    &       &   1   &   1   &   0   &   bias    &   0.079   &   0.079   &   -0.084  &   0.084   \\
    &       &       &       &       &   st.dev. &   (0.564) &   (0.615) &   (1.184) &   (0.991) \\
    &       &       &       &       &   mean s.e.   &   (0.543) &   (0.592) &   (1.155) &   (0.985) \\\cline{3-10}
    &       &   0   &   1   &   1   &   bias    &   0.038   &   0.111   &   -0.044  &   0.117   \\
    &       &       &       &       &   st.dev. &   (0.568) &   (0.601) &   (1.200) &   (0.995) \\
    &       &       &       &       &   mean s.e.   &   (0.547) &   (0.597) &   (1.132) &   (0.959) \\\hline
500 &   0.00    &   2   &   1   &   -1  &   bias    &   0.058   &   0.048   &   -0.176  &   0.165   \\
    &       &       &       &       &   st.dev. &   (0.407) &   (0.356) &   (0.897) &   (0.678) \\
    &       &       &       &       &   mean s.e.   &   (0.399) &   (0.350) &   (0.853) &   (0.662) \\\cline{3-10}
    &       &   1   &   1   &   0   &   bias    &   0.021   &   0.035   &   -0.004  &   0.018   \\
    &       &       &       &       &   st.dev. &   (0.354) &   (0.358) &   (0.787) &   (0.641) \\
    &       &       &       &       &   mean s.e.   &   (0.334) &   (0.348) &   (0.733) &   (0.609) \\\cline{3-10}
    &       &   0   &   1   &   1   &   bias    &   0.010   &   0.047   &   -0.004  &   0.041   \\
    &       &       &       &       &   st.dev. &   (0.343) &   (0.351) &   (0.733) &   (0.604) \\
    &       &       &       &       &   mean s.e.   &   (0.333) &   (0.349) &   (0.711) &   (0.590) \\\hline
500 &   0.75    &   2   &   1   &   -1  &   bias    &   0.055   &   0.011   &   -0.129  &   0.085   \\
    &       &       &       &       &   st.dev. &   (0.411) &   (0.366) &   (0.837) &   (0.647) \\
    &       &       &       &       &   mean s.e.   &   (0.387) &   (0.354) &   (0.799) &   (0.632) \\\cline{3-10}
    &       &   1   &   1   &   0   &   bias    &   0.007   &   0.022   &   0.022   &   -0.007  \\
    &       &       &       &       &   st.dev. &   (0.329) &   (0.360) &   (0.700) &   (0.575) \\
    &       &       &       &       &   mean s.e.   &   (0.327) &   (0.353) &   (0.693) &   (0.589) \\\cline{3-10}
    &       &   0   &   1   &   1   &   bias    &   0.019   &   0.023   &   -0.036  &   0.040   \\
    &       &       &       &       &   st.dev. &   (0.343) &   (0.370) &   (0.700) &   (0.581) \\
    &       &       &       &       &   mean s.e.   &   (0.333) &   (0.354) &   (0.686) &   (0.574) \\\hline\hline
\end{tabular}}
\caption{\em Simulation results for the proposed two-step estimator
based on 1000 samples of size $n=500,1000$ generated under different
models based on assumptions (\ref{simu1}) and (\ref{simu2}), with
$\rho=0.00,0.75$ and different values of $\al_1$ and $\be$ (in each
case $\al_2=1$ and $\de=0$).}\label{tab1}\vspace*{0.25cm}
\end{table}

\begin{table}[ht!]\centering
{\small\begin{tabular}{ccccclcccc}
\hline\hline        \vspace*{-0.4cm}                    \\
$n$ & $\rho$  &   $\al_1$   &   $\al_2$ &   $\be$ & & $\hat{\al}_1$
&   $\hat{\al}_2$   &   $\hat{\be}$   & $\hat{\de}$ \\
\hline
200 &   0.00    &   2   &   2   &   -1  &   bias    &   0.049   &   0.077   &   -0.075  &   0.103   \\
    &       &       &       &       &   st.dev. &   (0.593) &   (0.620) &   (1.363) &   (1.135) \\
    &       &       &       &       &   mean s.e.   &   (0.648) &   (0.726) &   (1.403) &   (1.193) \\\cline{3-10}
    &       &   1   &   2   &   0   &   bias    &   0.041   &   0.096   &   0.032   &   0.022   \\
    &       &       &       &       &   st.dev. &   (0.556) &   (0.646) &   (1.237) &   (1.066) \\
    &       &       &       &       &   mean s.e.   &   (0.550) &   (0.734) &   (1.212) &   (1.119) \\\cline{3-10}
    &       &   0   &   2   &   1   &   bias    &   0.053   &   0.140   &   -0.052  &   0.140   \\
    &       &       &       &       &   st.dev. &   (0.556) &   (0.627) &   (1.190) &   (1.020) \\
    &       &       &       &       &   mean s.e.   &   (0.553) &   (0.741) &   (1.190) &   (1.097) \\\hline
200 &   0.75    &   2   &   2   &   -1  &   bias    &   0.101   &   0.111   &   -0.165  &   0.175   \\
    &       &       &       &       &   st.dev. &   (0.610) &   (0.657) &   (1.330) &   (1.127) \\
    &       &       &       &       &   mean s.e.   &   (0.645) &   (0.745) &   (1.348) &   (1.171) \\\cline{3-10}
    &       &   1   &   2   &   0   &   bias    &   0.064   &   0.112   &   -0.016  &   0.064   \\
    &       &       &       &       &   st.dev. &   (0.591) &   (0.645) &   (1.225) &   (1.044) \\
    &       &       &       &       &   mean s.e.   &   (0.548) &   (0.743) &   (1.168) &   (1.096) \\\cline{3-10}
    &       &   0   &   2   &   1   &   bias    &   0.024   &   0.124   &   -0.013  &   0.114   \\
    &       &       &       &       &   st.dev. &   (0.575) &   (0.621) &   (1.162) &   (0.985) \\
    &       &       &       &       &   mean s.e.   &   (0.554) &   (0.749) &   (1.143) &   (1.068) \\\hline
500 &   0.00    &   2   &   2   &   -1  &   bias    &   0.065   &   0.060   &   -0.175  &   0.170   \\
    &       &       &       &       &   st.dev. &   (0.434) &   (0.476) &   (0.904) &   (0.735) \\
    &       &       &       &       &   mean s.e.   &   (0.399) &   (0.445) &   (0.840) &   (0.710) \\\cline{3-10}
    &       &   1   &   2   &   0   &   bias    &   0.030   &   0.077   &   -0.010  &   0.057   \\
    &       &       &       &       &   st.dev. &   (0.327) &   (0.445) &   (0.720) &   (0.671) \\
    &       &       &       &       &   mean s.e.   &   (0.334) &   (0.447) &   (0.731) &   (0.673) \\\cline{3-10}
    &       &   0   &   2   &   1   &   bias    &   0.024   &   0.084   &   -0.015  &   0.075   \\
    &       &       &       &       &   st.dev. &   (0.342) &   (0.487) &   (0.742) &   (0.697) \\
    &       &       &       &       &   mean s.e.   &   (0.333) &   (0.451) &   (0.716) &   (0.662) \\\hline
500 &   0.75    &   2   &   2   &   -1  &   bias    &   0.054   &   0.062   &   -0.134  &   0.142   \\
    &       &       &       &       &   st.dev. &   (0.413) &   (0.458) &   (0.832) &   (0.716) \\
    &       &       &       &       &   mean s.e.   &   (0.387) &   (0.451) &   (0.798) &   (0.693) \\\cline{3-10}
    &       &   1   &   2   &   0   &   bias    &   0.024   &   0.067   &   -0.012  &   0.056   \\
    &       &       &       &       &   st.dev. &   (0.350) &   (0.474) &   (0.723) &   (0.685) \\
    &       &       &       &       &   mean s.e.   &   (0.330) &   (0.451) &   (0.698) &   (0.656) \\\cline{3-10}
    &       &   0   &   2   &   1   &   bias    &   0.041   &   0.066   &   -0.064  &   0.088   \\
    &       &       &       &       &   st.dev. &   (0.336) &   (0.475) &   (0.701) &   (0.664) \\
    &       &       &       &       &   mean s.e.   &   (0.333) &   (0.451) &   (0.688) &   (0.642) \\\hline\hline
\end{tabular}}
\caption{\em Simulation results for the proposed two-step estimator
based on 1000 samples of size $n=500,1000$ generated under different
models based on assumptions (\ref{simu1}) and (\ref{simu2}), with
$\rho=0.00,0.75$ and different values of $\al_1$ and $\be$ (in each
case $\al_2=2$ and $\de=1$).}\label{tab2} \vspace*{0.1cm}
\end{table}

The simulations show that the estimators always have a very low bias
which, as may expected, tends to be smaller when the true value of
$\be$ is equal to 0 and when $n=500$ instead of $n=200$. It is also
worth noting that this bias is not considerably affected by $\rho$.
These conclusions are in agreement with those regarding the
asymptotic behavior of the estimator drawn on the basis of Figure
\ref{fig2}. Consequently, the rule of thumb that we use when all the
possible discordant configurations are not present seems to work
properly.

For what concerns the variability of the estimators, we observe that
the standard error of each of them is roughly proportional to
$1/\sqrt{n}$. This property is also in agreement with the asymptotic
results illustrated in Section \ref{sec:corr_est_asy}. Finally, for
each estimator, the average standard error is close enough to its
standard deviation. Relevant differences are only observed for
$n=200$ when the standard error occasionally tends to be larger than
the standard deviation. This confirms the validity of the proposed
estimator of the variance-covariance matrix of $\hat{\b\th}$ which
is based on the sandwich formula given in (\ref{sandwich}).
\section{Dealing with missing responses}
We now illustrate how the proposed causal model and the two-step
conditional logistic estimator may be extended to the case of
missing responses. For this aim, we introduce the binary indicator
$R_h$, $h=1,2$, equal to 1 if the response variable $Y_h$ is
observable and to 0 otherwise.

\subsection{Causal model}

With missing responses, we extend the model introduced in Section
\ref{sec:ass} by assuming that:
\begin{itemize}
\item[B1:] $R_1\ind Y_1|(U,\b V)$;
\item[B2:] $C\ind(R_1,Y_1)|(U,\b V)$;
\item[B3:] $Z\ind(U,Y_1,R_1,C)|\b V$;
\item[B4:] $X\ind(U,\b V,Y_1,R_1)|(C,Z)$ and, with probability 1,
$X=Z$ when $C=1$ (compliers) and $X=0$ when $C=0$ (never-takers);
\item[B5:] $Y_2\ind (Y_1,R_1,Z)|(U,\b V,C,X)$;
\item[B6:] $R_2\ind(Y_1,R_1,Z,Y_2)|(U,\b V,C,X)$;
\item[B7:] for any $u,\b v,c$ and $x$, we have
\[
\logit[p(Y_2=1|u,\b v,c,x)]-\logit[p(Y_1=1|u,\b v)]=\b a(\b
v,x)\tr\b\al+c(1-x)\b b(\b v)\tr\b\be,
\]
with the functions $\b a(\b v,x)$ and $\b b(\b v)$ and the
corresponding parameter vectors defined as in Section 2.1.
\end{itemize}

New assumptions are essentially B1 and B6 concerning the conditional
independence between $R_1$ and $Y_1$ and that between $R_2$ and
$Y_2$ given observable and unobservable covariates. This assumption
is weaker than the assumption that responses are missing at random
given the observable covariates \cite{Rubin:1976}.

\begin{figure}[ht]\centering
\epsfig{figure=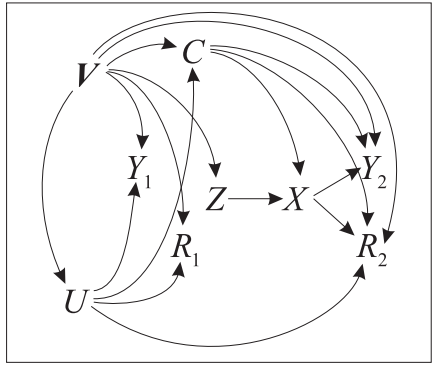,height=4.8439cm} \caption{\em DAG for
the model based on assumptions B1-B7. $U$ and $\b V$ represent
unobservable and observable covariates affecting the response
variables $Y_1$ and $Y_2$ and the indicator variables $R_1$ and
$R_2$ for the response variables being observable, $C$ is a binary
variable for the compliance status and $Z$ and $X$ are binary
variables for the assigned and the received
treatment.}\label{fig4}\vspace*{0.25cm}
\end{figure}

The resulting model is represented by the DAG in Figure \ref{fig4}.

\subsection{Two-step conditional logistic estimator}

Under assumptions B1-B7, we have that
\[
p(y_1,r_1,z,x,y_2,r_2|u,\b v,c)=p(y_1|u,\b v)p(r_1|u,\b v)p(z|\b
v)p(x|c,z)p(y_2|u,\b v,c,x)p(r_2|u,\b v,c,x).
\]
Then, marginalizing with respect to $C$ and $U$, we find that the
probability mass function of the conditional distribution of
$(Y_1,R_1,Z,X,Y_2,R_2)$ given $\b V$ is equal to
\begin{eqnarray*}
f^*(y_1,r_1,z,x,y_2,r_2|\b v)&=&p(z|\b v)\int\frac{e^{y_+\la(u,\bl v)}}{1+e^{\la(u,\bl v)}}p(r_1|u,\b v)\times\nonumber\\
&\times&\sum_c p(x|c,z)\frac{e^{y_2t(c,\bl v,x)}}{1+e^{\la(u,\bl
v)+t(c,\bl v,x)}}p(r_2|u,\b v,c,x)\pi(c|u,\b v)\phi(u|\b v)du,
\end{eqnarray*}
where, as in Section \ref{sec:cond}, $t(c,\b v,x)=\b a(\b
v,x)\tr\b\al+c(1-x)\b b(\b v)\tr\b\be$.

When $z=1$, the above expression simplifies to
\begin{eqnarray*}
f^*(y_1,r_1,1,x,y_2,r_2|\b v)&=&p(Z=1|\b
v)\times\\
&\times&\int\frac{e^{y_+\la(u,\bl v)}}{1+e^{\la(u,\bl v)}}p(r_1|u,\b
v)\frac{e^{y_2\bl a(\bl v,x)\tr\bl\al}}{1+e^{\la(u,\bl v)+\bl a(\bl
v,x)\tr\bl\al}}p(r_2|u,\b v,c,x)\phi(u|\b v)du,
\end{eqnarray*}
and this implies that
\[
f^*(y_2|\b v,1,x)=\frac{f^*(1-y_2,1,1,x,y_2,1|\b
v)}{f^*(0,1,1,x,1,1|\b v)+f^*(1,1,1,x,0,1|\b v)}=\frac{e^{y_2\bl
a(\bl v,x)\tr\bl\al}}{1+e^{\bl a(\bl v,x)\tr\bl\al}},
\]
where, in general, $f^*(y_2|\b v,z,x)=p(y_2|\b
v,R_1=1,z,x,R_2=1,Y_+=1)$. This is due to the definition of
$p(x|c,z)$ which, as already noted, may only be equal to 0 or 1. The
same does not happen when $z=0$ since in this case
\begin{eqnarray*}
f^*(y_1,r_1,0,0,y_2,r_2|\b v)&=&p(Z=0|\b v)g^*(y_1,r_1,y_2,r_2|b v)\\
g^*(y_1,r_1,y_2,r_2|\b v)&=&\int\frac{e^{y_+\la(u,\bl
v)}}{1+e^{\la(u,\bl v)}}p(r_1|u,\b v)\times\\
&\times&\sum_c \frac{e^{y_2t(c,\bl v,0)}}{1+e^{\la(u,\bl v)+t(c,\bl
v,0)}}p(r_2|u,\b v,c,0)\pi(c|u,\b v)\phi(u|\b v)du.
\end{eqnarray*}
Then, as in Section \ref{sec:corr}, we consider a first-order Taylor
series expansion of $\log g^*(y_1,r_1,y_2,r_2|\b v)$ around
$\b\be=\b 0$ and we find that
\[
\log g^*(y_1,r_1,y_2,r_2|\b v)\approx\log g^*_0(y_1,r_1,y_2,r_2|\b
v)+ h^*(y_1,r_1,y_2,r_2|\b v)\b b(\b v)\tr\b\be
\]
where $g_0^*(y_1,r_1,y_2,r_2|\b v)$ is the function
$g^*(y_1,r_1,y_2,r_2|\b v)$ computed at $\b\be=\b 0$ and
\begin{eqnarray*}
h^*(y_1,r_1,y_2,r_2|\b v)&=&\frac{2y_2-1}{g_0^*(y_1,r_1,y_2,r_2|\b
v)}\int\frac{e^{y_+\la(u,\bl v)}}{1+e^{\la(u,\bl v)}}p(r_1|u,\b
v)\times\\&\times&\frac{e^{t(0,\bl v,0)}}{[1+e^{\la(u,\bl v)+t(0,\bl
v,0)}]^2}p(r_2|u,\b v,c,0)\pi(1|u,\b v)\phi(u|\b v)du.
\end{eqnarray*}
Consequently, we have
\[
\log\frac{f^*(0,1,0,0,1,1|\b v)}{f^*(1,1,0,0,0,1|\b v)}\approx\b
a(\b v,0)\tr\b\al+h^*(\b v)\b b(\b v)\tr\b\be,
\]
where $h^*(\b v)$ is a correction factor which is simply equal to
$\pi(1|\b v)=p(C=1|\b v)$ when $C$ is conditionally independent of
$U$ given $\b V$. Finally, we have
\[
f^*(y_2|\b v,0,0)\approx\frac{e^{\bl a(\bl v,0)\tr\bl\al+h(\bl v)\bl
b(\bl v)\tr\bl\be}}{1+e^{\bl a(\bl v,0)\tr\bl\al+h(\bl v)\bl b(\bl
v)\tr\bl\be}},
\]
which is exactly the same expression given in (\ref{eq:approx}).

On the basis of the above arguments, in the case of missing
responses we propose to use a two-step estimator which has the same
structure as that described in Section \ref{sec:corr_est} and is
based on a logit model of type (\ref{logit}) for the probability
$\pi(1|\b v)$ of being a complier given the covariates. Here, for
each subject $i$, with $i=1,\ldots,n$, we observe $\b v_i$,
$r_{i1}$, $z_i$, $x_i$ and $r_{i2}$, where, for $h=1,2$, $r_{hi}$ is
the observed value of $R_h$. For the same subject we also observe
$y_{i1}$ if $r_{i1}=1$ and $y_{i2}$ if $r_{i2}=1$.

The estimator is based on the following steps which, as for the
initial estimator, may be performed on the basis of standard
estimation algorithms:

\begin{enumerate}
\item {\em Estimation of $\b\eta$}. This is based on the observed
values $\b v_i$ and $x_i$ for every $i$ such that $z_i=1$ and
proceeds by maximizing the weighted log-likelihood
\[
\ell_1^*(\b\eta)=\sum_i\frac{z_i}{p(z_i|\b v_i)}[x_i\log\pi(1|\b
v_i)+(1-x_i)\log\pi(0|\b v_i)].
\]
\item {\em Estimation of $\b\al$ and $\b\be$}. This is performed by maximizing
the conditional log-likelihood
\[
\ell_2(\b\al,\b\be;\b\eta) = \sum_i
d_ir_{i1}r_{i2}\ell_{i2}(\b\al,\b\be;\b\eta),
\]
where  $\ell_{i2}(\b\al,\b\be;\b\eta)$ is defined as in (\ref{lk2}).
\end{enumerate}

Note that the only difference with respect to the estimator in
Section \ref{sec:corr_est} is in the second step where we consider
only the subjects who respond at both occasions, whereas at the
first step we consider all subjects in order to estimate the model
for the probability of being a complier. Moreover, standard errors
for the estimator can again be computed on the basis of the sandwich
formula (\ref{sandwich}).
\subsection{Properties of the two-step estimator}
A final point concerns the properties of the estimators
$\hat{\b\al}$, $\hat{\b\be}$ and $\hat{\de}(\b v)$ with missing
responses. These estimators have the same asymptotic properties they
have when the response variables are always observable (see Section
\ref{sec:corr_est_asy}). The main result is that these estimators
are consistent when $\b\be_0=\b 0$, regardless of the
parametrization used in the logit model (\ref{logit}) for the
probability to comply. When $\b\be_0\neq\b 0$, the estimators
$\hat{\b\al}$ and $\hat{\b\be}$ converge to $\b\al_*$ and $\b\be_*$,
respectively, as $n\rightarrow\infty$. These limits are equal to the
true values $\b\al_0$ and $\b\be_0$ when $\b\be_0=\b 0$ and are
expected to be reasonably close to these true values when $\b\be_0$
is not too far from $\b 0$. The same may be said for the estimator
$\hat{\de}(\b v)$ which converges to $\de_*(\b v)$.

To illustrate the above point, in Figure 3 we report some plots of
$\de_*$ with respect to $\be$ under a true model involving only one
observable covariate and based on the same assumptions illustrated
in Section \ref{sec:corr_est_asy}, see in particular (\ref{simu1})
and (\ref{simu2}), beyond the assumption that $R_1$ and $R_2$ have
Bernoulli distribution with probabilities of success chosen,
respectively, as follows
\begin{equation}
\begin{array}{rcl}
p(R_1=1|u,v) &=& \expit[1+(u+v)/\sqrt{1+\rho}/2],\vspace*{0.25cm}\\
p(R_2=1|u,v,c,x) &=& \expit[1+(u+v)/\sqrt{1+\rho}/2+c/2+x/2].
\end{array}
\label{simu4}
\end{equation}
The estimators we considered are:
\begin{itemize}
\item $\hat{\de}_{null}$: two-step conditional logistic
estimator of $\de$ based on a model for the probability to comply of
type (\ref{logit}) with $m(v)=1$;
\item $\hat{\de}_{cov}$: as above with
$\b m(v)=(1,v)\tr$, so that the covariate is also used to predict
the probability to comply;
\item $\hat{\de}_{itt}$: ITT estimator based on the conditional logistic regression on only the
subjects who respond at both occasions, i.e. we regress $Y_2$ on
$(1,Z)$ given $Y_+=1$, $R_1=1$ and $R_2=1$;
\item $\hat{\de}_{tr}$: TR estimator based the conditional logistic
regression of $Y_2$ on $(1,X)$ given $Y_+=1$, $R_1=1$ and $R_2=1$.
\end{itemize}
The resulting plots closely resemble those in Figure 2 and then
similar conclusions may be drawn about the proposed estimator. In
particular, we again note the small distance between the limit in
probability of the estimator and the true value of the parameter.

\begin{figure}[ht]\centering
\epsfig{figure=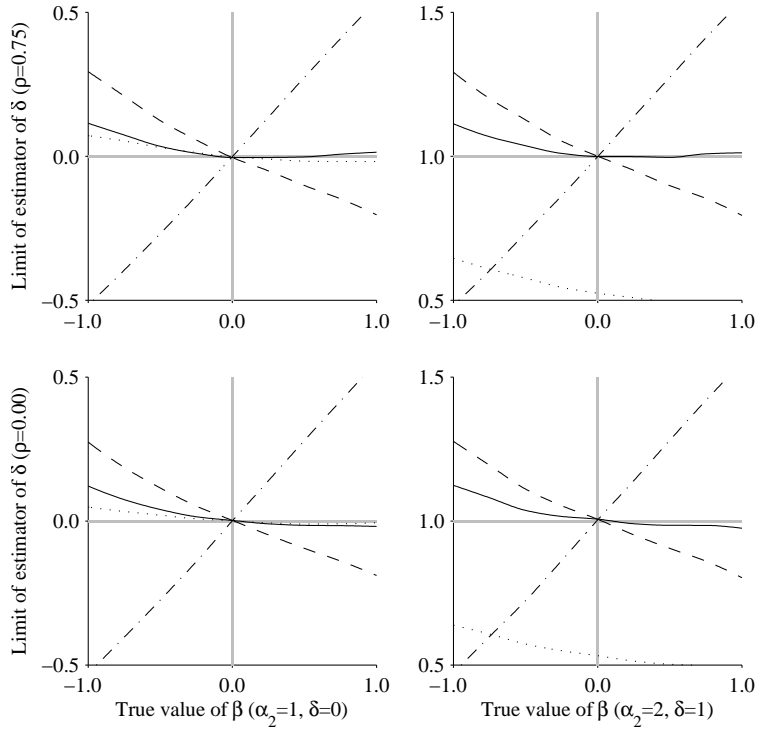,width=10cm} \caption{\em Limit
of the estimators $\hat{\de}_{null}$ (dashed line),
$\hat{\de}_{cov}$ (solid line), $\hat{\de}_{itt}$ (dotted line) and
$\hat{\de}_{tr}$ (dash-dotted line) under assumptions (\ref{simu1}),
(\ref{simu2}) and (\ref{simu4}), with $\rho=0.00,0.75$, $\be$
between -1 and 1 and $\de=0,1$, with $\al_1$ and $\al_2$ defined
accordingly.}\vspace*{0.25cm}
\end{figure}

Under the same true model assumed above, we studied by simulation
the finite-sample properties of the estimators $\hat{\b\al}$,
$\hat{\b\be}$ and $\hat{\de}$. As usual, we focused on the
estimators which exploit the covariate to predict the probability to
comply, and then we let $\b m(\b v)=(1,v)\tr$ in (\ref{logit}).
Under the same setting of the simulations in Section 4.2.2, we
obtained the results reported in Tables \ref{tab3} and \ref{tab4}
when (\ref{simu4}) is assumed. These results are very similar to
those reported in Tables \ref{tab1} and \ref{tab2} for the case in
which the response variables are always observed. The main
difference is in the variability of the estimators which is
obviously larger because of the presence of missing responses.

\begin{table}[ht!]\centering
{\small\begin{tabular}{ccccclcccc}
\hline\hline        \vspace*{-0.4cm}                    \\
$n$ & $\rho$  &   $\al_1$   &   $\al_2$ &   $\be$ & & $\hat{\al}_1$
&   $\hat{\al}_2$   &   $\hat{\be}$   & $\hat{\de}$ \\
\hline
200 &   0.00    &   2   &   1   &   -1  &   bias    &   -0.107  &   0.076   &   0.208   &   -0.025  \\
    &       &       &       &       &   st.dev. &   (0.628) &   (0.706) &   (1.573) &   (1.404) \\
    &       &       &       &       &   mean s.e.   &   (0.910) &   (0.733) &   (1.910) &   (1.445) \\\cline{3-10}
    &       &   1   &   1   &   0   &   bias    &   0.053   &   0.129   &   0.024   &   0.053   \\
    &       &       &       &       &   st.dev. &   (0.710) &   (0.723) &   (1.624) &   (1.369) \\
    &       &       &       &       &   mean s.e.   &   (0.825) &   (0.742) &   (1.736) &   (1.356) \\\cline{3-10}
    &       &   0   &   1   &   1   &   bias    &   0.045   &   0.071   &   0.061   &   -0.034  \\
    &       &       &       &       &   st.dev. &   (0.776) &   (0.683) &   (1.654) &   (1.310) \\
    &       &       &       &       &   mean s.e.   &   (0.813) &   (0.725) &   (1.686) &   (1.303) \\\hline
    200 &   0.75    &   2   &   1   &   -1  &   bias    &   -0.069  &   0.065   &   0.161   &   -0.028  \\
    &       &       &       &       &   st.dev. &   (0.639) &   (0.696) &   (1.457) &   (1.310) \\
    &       &       &       &       &   mean s.e.   &   (0.928) &   (0.752) &   (1.841) &   (1.381) \\\cline{3-10}
    &       &   1   &   1   &   0   &   bias    &   0.078   &   0.069   &   -0.035  &   0.026   \\
    &       &       &       &       &   st.dev. &   (0.720) &   (0.699) &   (1.515) &   (1.213) \\
    &       &       &       &       &   mean s.e.   &   (0.813) &   (0.751) &   (1.625) &   (1.283) \\\cline{3-10}
    &       &   0   &   1   &   1   &   bias    &   0.003   &   0.076   &   0.067   &   0.007   \\
    &       &       &       &       &   st.dev. &   (0.801) &   (0.696) &   (1.613) &   (1.256) \\
    &       &       &       &       &   mean s.e.   &   (0.816) &   (0.750) &   (1.575) &   (1.228) \\\hline
500 &   0.00    &   2   &   1   &   -1  &   bias    &   0.113   &   0.070   &   -0.272  &   0.230   \\
    &       &       &       &       &   st.dev. &   (0.594) &   (0.467) &   (1.235) &   (0.907) \\
    &       &       &       &       &   mean s.e.   &   (0.608) &   (0.440) &   (1.211) &   (0.862) \\\cline{3-10}
    &       &   1   &   1   &   0   &   bias    &   0.023   &   0.038   &   -0.008  &   0.023   \\
    &       &       &       &       &   st.dev. &   (0.512) &   (0.457) &   (1.022) &   (0.809) \\
    &       &       &       &       &   mean s.e.   &   (0.482) &   (0.438) &   (0.998) &   (0.781) \\\cline{3-10}
    &       &   0   &   1   &   1   &   bias    &   0.036   &   0.008   &   0.028   &   -0.056  \\
    &       &       &       &       &   st.dev. &   (0.482) &   (0.443) &   (1.009) &   (0.785) \\
    &       &       &       &       &   mean s.e.   &   (0.476) &   (0.433) &   (0.976) &   (0.759) \\\hline
500 &   0.75    &   2   &   1   &   -1  &   bias    &   0.088   &   0.037   &   -0.187  &   0.135   \\
    &       &       &       &       &   st.dev. &   (0.578) &   (0.463) &   (1.123) &   (0.844) \\
    &       &       &       &       &   mean s.e.   &   (0.596) &   (0.451) &   (1.142) &   (0.823) \\\cline{3-10}
    &       &   1   &   1   &   0   &   bias    &   0.042   &   0.064   &   -0.025  &   0.047   \\
    &       &       &       &       &   st.dev. &   (0.530) &   (0.468) &   (1.045) &   (0.799) \\
    &       &       &       &       &   mean s.e.   &   (0.490) &   (0.454) &   (0.964) &   (0.758) \\\cline{3-10}
    &       &   0   &   1   &   1   &   bias    &   0.035   &   0.035   &   0.001   &   -0.001  \\
    &       &       &       &       &   st.dev. &   (0.503) &   (0.482) &   (0.976) &   (0.769) \\
    &       &       &       &       &   mean s.e.   &   (0.483) &   (0.454) &   (0.937) &   (0.735) \\\hline\hline
\end{tabular}}
\caption{\em Simulation results for the proposed two-step estimator
based on 1000 samples of size $n=500,1000$ generated under different
models based on assumptions (\ref{simu1}) and (\ref{simu2}) and
(\ref{simu4}), with $\rho=0.00,0.75$ and different values of $\al_1$
and $\be$ (in each case $\al_2=1$ and
$\de=0$).}\label{tab3}\vspace*{0.25cm}
\end{table}

\begin{table}[ht!]\centering
{\small\begin{tabular}{ccccclcccc}
\hline\hline        \vspace*{-0.4cm}                    \\
$n$ & $\rho$  &   $\al_1$   &   $\al_2$ &   $\be$ & & $\hat{\al}_1$
&   $\hat{\al}_2$   &   $\hat{\be}$   & $\hat{\de}$ \\
\hline
200 &   0.00    &   2   &   2   &   -1  &   bias    &   -0.084  &   -0.019  &   0.227   &   -0.162  \\
    &       &       &       &       &   st.dev. &   (0.590) &   (0.614) &   (1.547) &   (1.329) \\
    &       &       &       &       &   mean s.e.   &   (0.913) &   (0.876) &   (1.915) &   (1.531) \\\cline{3-10}
    &       &   1   &   2   &   0   &   bias    &   0.056   &   0.010   &   -0.046  &   0.000   \\
    &       &       &       &       &   st.dev. &   (0.733) &   (0.614) &   (1.644) &   (1.307) \\
    &       &       &       &       &   mean s.e.   &   (0.818) &   (0.879) &   (1.713) &   (1.426) \\\cline{3-10}
    &       &   0   &   2   &   1   &   bias    &   0.043   &   -0.021  &   0.070   &   -0.134  \\
    &       &       &       &       &   st.dev. &   (0.784) &   (0.589) &   (1.650) &   (1.236) \\
    &       &       &       &       &   mean s.e.   &   (0.815) &   (0.871) &   (1.686) &   (1.389) \\\hline
200 &   0.75    &   2   &   2   &   -1  &   bias    &   -0.092  &   -0.022  &   0.188   &   -0.118  \\
    &       &       &       &       &   st.dev. &   (0.617) &   (0.604) &   (1.436) &   (1.239) \\
    &       &       &       &       &   mean s.e.   &   (0.912) &   (0.890) &   (1.827) &   (1.468) \\\cline{3-10}
    &       &   1   &   2   &   0   &   bias    &   0.033   &   -0.014  &   0.101   &   -0.148  \\
    &       &       &       &       &   st.dev. &   (0.714) &   (0.607) &   (1.513) &   (1.198) \\
    &       &       &       &       &   mean s.e.   &   (0.799) &   (0.892) &   (1.595) &   (1.364) \\\cline{3-10}
    &       &   0   &   2   &   1   &   bias    &   0.031   &   -0.037  &   0.021   &   -0.089  \\
    &       &       &       &       &   st.dev. &   (0.801) &   (0.591) &   (1.549) &   (1.137) \\
    &       &       &       &       &   mean s.e.   &   (0.819) &   (0.880) &   (1.593) &   (1.322) \\\hline
500 &   0.00    &   2   &   2   &   -1  &   bias    &   0.071   &   0.129   &   -0.186  &   0.243   \\
    &       &       &       &       &   st.dev. &   (0.601) &   (0.602) &   (1.256) &   (1.005) \\
    &       &       &       &       &   mean s.e.   &   (0.599) &   (0.587) &   (1.202) &   (0.954) \\\cline{3-10}
    &       &   1   &   2   &   0   &   bias    &   0.065   &   0.115   &   -0.073  &   0.122   \\
    &       &       &       &       &   st.dev. &   (0.505) &   (0.598) &   (1.029) &   (0.883) \\
    &       &       &       &       &   mean s.e.   &   (0.490) &   (0.582) &   (1.012) &   (0.881) \\\cline{3-10}
    &       &   0   &   2   &   1   &   bias    &   0.023   &   0.138   &   0.049   &   0.066   \\
    &       &       &       &       &   st.dev. &   (0.494) &   (0.583) &   (0.995) &   (0.875) \\
    &       &       &       &       &   mean s.e.   &   (0.475) &   (0.586) &   (0.978) &   (0.863) \\\hline
500 &   0.75    &   2   &   2   &   -1  &   bias    &   0.112   &   0.123   &   -0.263  &   0.274   \\
    &       &       &       &       &   st.dev. &   (0.583) &   (0.586) &   (1.148) &   (0.933) \\
    &       &       &       &       &   mean s.e.   &   (0.598) &   (0.594) &   (1.140) &   (0.912) \\\cline{3-10}
    &       &   1   &   2   &   0   &   bias    &   0.034   &   0.098   &   -0.033  &   0.097   \\
    &       &       &       &       &   st.dev. &   (0.507) &   (0.588) &   (1.009) &   (0.848) \\
    &       &       &       &       &   mean s.e.   &   (0.484) &   (0.591) &   (0.952) &   (0.849) \\\cline{3-10}
    &       &   0   &   2   &   1   &   bias    &   0.060   &   0.130   &   -0.085  &   0.155   \\
    &       &       &       &       &   st.dev. &   (0.476) &   (0.612) &   (0.917) &   (0.827) \\
    &       &       &       &       &   mean s.e.   &   (0.477) &   (0.601) &   (0.920) &   (0.832) \\\hline\hline
\end{tabular}}
\caption{\em Simulation results for the proposed two-step estimator
based on 1000 samples of size $n=500,1000$ generated under different
models based on assumptions (\ref{simu1}) and (\ref{simu2}) and
(\ref{simu4}), with $\rho=0.00,0.75$ and different values of $\al_1$
and $\be$ (in each case $\al_2=2$ and
$\de=1$).}\label{tab4}\vspace*{0.25cm}
\end{table}

\newpage\section{An application}
To illustrate the approach proposed in this paper, we analyzed the
dataset coming from the randomized experiment on BSE already
mentioned in Section 1.

The study took place between the beginning of 1988 and the end of
1990 at the Oncologic Center of the Faenza District, Italy. The
sample used in the study consists of 657 women aged 20 to 64 years,
who were randomly assigned to the control, consisting of learning
how to perform BSE through a standard method, or to the treatment,
consisting of a training course held by a specialized medical staff.
Only women assigned to the treatment can access it and then
non-compliance may be observed only among these subjects. In
particular, of the 330 women randomly assigned to the treatment, 182
attended the course and so they may be considered as compliers; the
remaining women may be considered as never-takes.\newpage

The efficacy of the treatment is measured by two binary response
variables, observed before and after the treatment/control, which
indicate if {\em BSE is regularly practised} and if the {\em quality
of BSE practise is adequate}. Several covariates are also available,
such as {\em age}, {\em number of children}, {\em educational
level}, {\em occupational status}, {\em presence of previous cancer
pathologies in the woman or her family}, {\em menopause} and {\em
adequate knowledge of breast pathophysiology}. Finally, some
response variables are not observed and these have to be treated as
missing.

The dataset has already been analyzed by \cite{Ferro:1996}, on the
basis of a standard conditional logistic approach, and by
\cite{Mealli:2004}, who exploited a potential outcome approach
allowing for missing responses, which is related to that of
\cite{Frang:1999}.

In analyzing the dataset, we first considered the effect of the
treatment on practicing BSE. In this case, $Y_1$ is equal to 1 if a
woman regularly practises BSE before the treatment and to 0
otherwise. Similarly, $Y_2$ is equal to 1 if a woman regularly
practises BSE after the treatment and to 0 otherwise. The first
variable was observed for the 93.61\% of the sample and the second
for the 65.30\%. We then followed the method for missing responses
described in Section 5 under assumption (\ref{par2}) for the
parametrization of the causal model. In particular, we first
computed the estimators $\hat{\b\al}_{null}$, $\hat{\be}_{null}$ and
$\hat{\de}_{null}$, based on predicting the probability to comply
only on the basis of the indicator variable for the second response
variable being observable, and the estimators $\hat{\b\al}_{cov}$,
$\hat{\be}_{cov}$ and $\hat{\de}_{cov}$, which also consider the
covariates age and age-squared in the model used to predict this
probability. These covariates are included since are among those
with the most significant effect on the probability to comply. We
also considered the ITT estimator $\hat{\de}_{itt}$ and the TR
estimator $\hat{\de}_{tr}$ defined as in Section 5.3. The results
are displayed in Table \ref{tab5}.

Our first conclusion is that the inclusion of the covariates in the
model for the probability to comply does not dramatically affect the
estimates of the parameters and of the causal effect computed
following our approach. In particular, the estimate of $\al_2$
remains unchanged by the inclusion of these covariates, since this
estimate exploits only the data deriving from the treatment arm.
Overall, we can observe an effect of the control on never-takers,
corresponding to $\al_1$, which is not significant. Moreover, the
estimate of the parameter $\be$ is very different from zero,
indicating a great difference between compliers and never-takers for
what concerns this effect. Then, we conclude that the effect of the
treatment over control on practicing BSE ($\de$) is not significant.
A similar conclusion is reached on the basis of the ITT estimator of
$\de$, whereas the TR estimator attains a value much higher of that
of the other estimators, since it does not distinguish between
compliers and never-takers for what concerns the effect of the
treatment.\newpage

\begin{table}[ht]\centering
{\small\begin{tabular}{rrrrrrrrrr}\hline\hline\vspace*{-0.4cm}\\
\multicolumn4c{$\hat{\b\al}_{null}$} &&
\multicolumn4c{$\hat{\b\al}_{cov}$}\\\cline{1-4}\cline{6-9}
\multicolumn1c{est.} & \multicolumn1c{se} & \multicolumn1c{$t$} &
\multicolumn1c{$p$-value} &&
\multicolumn1c{est.} & \multicolumn1c{se} & \multicolumn1c{$t$} & \multicolumn1c{$p$-value}\\
\cline{1-4}\cline{6-9}
   -0.336 & 0.586 & -0.575 & 0.566 &&  -0.143 & 0.573 & -0.250 & 0.802\\
    2.241 & 0.470 &  4.763 & 0.000 &&   2.241 & 0.471 &  4.763 & 0.000\\
 \hline\hline\vspace*{-0.4cm}\\
\multicolumn4c{$\hat{\be}_{null}$} &&
\multicolumn4c{$\hat{\be}_{cov}$}\\\cline{1-4}\cline{6-9}
\multicolumn1c{est.} & \multicolumn1c{se} & \multicolumn1c{$t$} &
\multicolumn1c{$p$-value} &&
\multicolumn1c{est.} & \multicolumn1c{se} & \multicolumn1c{$t$} & \multicolumn1c{$p$-value}\\
\cline{1-4}\cline{6-9}
       2.779 & 1.207 & 2.302 & 0.021 &&  2.382 & 1.244 & 1.914 & 0.056\\\hline\hline\vspace*{-0.4cm}\\
\multicolumn4c{$\hat{\de}_{null}$} &&
\multicolumn4c{$\hat{\de}_{cov}$}\\\cline{1-4}\cline{6-9}
\multicolumn1c{est.} & \multicolumn1c{se} & \multicolumn1c{$t$} &
\multicolumn1c{$p$-value} &&
\multicolumn1c{est.} & \multicolumn1c{se} & \multicolumn1c{$t$} & \multicolumn1c{$p$-value}\\
\cline{1-4}\cline{6-9}
  -0.202 & 0.882 & -0.229 & 0.819 &&   0.002 & 0.909 & 0.003 & 0.998
  \\\hline\hline\vspace*{-0.4cm}\\
\multicolumn4c{$\hat{\de}_{itt}$} &&
\multicolumn4c{$\hat{\de}_{tr}$}\\\cline{1-4}\cline{6-9}
\multicolumn1c{est.} & \multicolumn1c{se} & \multicolumn1c{$t$} &
\multicolumn1c{$p$-value} &&
\multicolumn1c{est.} & \multicolumn1c{se} & \multicolumn1c{$t$} & \multicolumn1c{$p$-value}\\
\cline{1-4}\cline{6-9} 0.2701 & 0.4503 & 0.5998 &
0.5486 && 1.3652 & 0.5405 & 2.5260 & 0.0115\\
\hline\hline
\end{tabular}}
\caption{\em Estimates of the causal parameters obtained on the
basis of the proposed approach and the ITT and TR approaches when
the response variable is equal to 1 for a woman regularly practicing
BSE and to 0 otherwise.}\label{tab5}\vspace*{0.25cm}
\end{table}

We then considered the effect of the treatment on the quality of the
BSE practise. As in \cite{Mealli:2002}, this is measured through the
binary response variables $Y_1$ and $Y_2$ (here redefined) which are
equal to 1 if the score assigned by the medical staff to the quality
of the BSE practise is greater than the sample median and to 0
otherwise. As usual, $Y_1$ is a pre-treatment variable and $Y_2$ is
a post-treatment variable. Obviously, these variables are observable
only if BSE is practised and so we again used the method for missing
responses described in Section 5. In particular, $Y_1$ was observed
for the 54.80\% of the sample and $Y_2$ for the 51.93\%. The results
obtained from the application of the same estimators mentioned above
are reported in Table \ref{tab6}.

In this case, the inclusion of the covariates age and age-squared in
the model for predicting the probability to comply has a slight
effect on the estimates of the parameters $\b\al$, $\be$ and $\de$
computed on the basis of the proposed approach. Never-takers and
compliers now appear to be less distant in terms of reaction to the
control, whose effect is not significant for both subpopulations. On
the other hand, the effect of the treatment on compliers is
significant as well as the causal effect of the treatment over
control. The estimate of this causal effect is in this case close to
the RT estimate, whereas the ITT estimate is much smaller, even if
it remains significantly greater than 0.

\begin{table}[ht]\centering
{\small\begin{tabular}{rrrrrrrrrr}\hline\hline\vspace*{-0.4cm}\\
\multicolumn4c{$\hat{\b\al}_{null}$} &&
\multicolumn4c{$\hat{\b\al}_{cov}$}\\\cline{1-4}\cline{6-9}
\multicolumn1c{est.} & \multicolumn1c{se} & \multicolumn1c{$t$} &
\multicolumn1c{$p$-value} &&
\multicolumn1c{est.} & \multicolumn1c{se} & \multicolumn1c{$t$} & \multicolumn1c{$p$-value}\\
\cline{1-4}\cline{6-9}
    0.000 & 0.707 & 0.000 & 1.000 &&  -0.112 & 0.687 & -0.164 & 0.870\\
    3.611 & 1.013 & 3.563 & 0.000 &&   3.611 & 1.013 &  3.563 & 0.000\\
\hline\hline\vspace*{-0.4cm}\\
\multicolumn4c{$\hat{\be}_{null}$} &&
\multicolumn4c{$\hat{\be}_{cov}$}\\\cline{1-4}\cline{6-9}
\multicolumn1c{est.} & \multicolumn1c{se} & \multicolumn1c{$t$} &
\multicolumn1c{$p$-value} &&
\multicolumn1c{est.} & \multicolumn1c{se} & \multicolumn1c{$t$} & \multicolumn1c{$p$-value}\\
\cline{1-4}\cline{6-9}
      0.880 & 1.406 & 0.626 & 0.531 && 1.103 & 1.323 & 0.834 & 0.404\\
\hline\hline\vspace*{-0.4cm}\\\multicolumn4c{$\hat{\de}_{null}$} &&
\multicolumn4c{$\hat{\de}_{cov}$}\\\cline{1-4}\cline{6-9}
\multicolumn1c{est.} & \multicolumn1c{se} & \multicolumn1c{$t$} &
\multicolumn1c{$p$-value} &&
\multicolumn1c{est.} & \multicolumn1c{se} & \multicolumn1c{$t$} & \multicolumn1c{$p$-value}\\
\cline{1-4}\cline{6-9}
      2.731 & 1.300 & 2.100 & 0.036 &&   2.620 & 1.262 & 2.076 & 0.038\\
\hline\hline\vspace*{-0.4cm}\\
\multicolumn4c{$\hat{\de}_{itt}$} &&
\multicolumn4c{$\hat{\de}_{tr}$}\\\cline{1-4}\cline{6-9}
\multicolumn1c{est.} & \multicolumn1c{se} & \multicolumn1c{$t$} &
\multicolumn1c{$p$-value} &&
\multicolumn1c{est.} & \multicolumn1c{se} & \multicolumn1c{$t$} & \multicolumn1c{$p$-value}\\
\cline{1-4}\cline{6-9}  1.6186 & 0.5704 & 2.8377 & 0.0045 &&  3.2055 & 1.0537 &  3.0420 & 0.0024\\
\hline\hline
\end{tabular}}
\caption{\em Estimates of the causal parameters obtained on the
basis of the proposed approach and the ITT and TR approaches when
the response variable is equal to 1 when the quality of the BSE
practise is adequate and to 0
otherwise.}\label{tab6}\vspace*{0.25cm}
\end{table}

Overall, the results obtained with the proposed approach are in
accordance with those of \cite{Ferro:1996} and \cite{Mealli:2004},
who concluded that the training course has not a significant effect
on practising BSE, but it has a significant effect on the quality of
the BSE practise.
\section{Discussion}\label{sec:conc}
A causal model has been introduced to study the behavior of the
conditional logistic estimator as a tool of analysis of data coming
from two-arm experimental studies with possible non-compliance. The
model is applicable with binary outcomes observed before and after
the treatment (or control). It is formulated on the basis of latent
variables for the effect of unobservable covariates at both
occasions and to account for the difference between compliers and
never-takers in terms of reaction to control and treatment. A
correction for the bias of the conditional logistic estimator has
also been proposed which can be exploited when we want to estimate
the causal effect of the treatment over control in the subpopulation
of compliers. It results a two-step estimator which has some
connection with the estimators of \cite{Heckman:1979} and
\cite{Nag:2000} and represents an extension of the standard
conditional logistic estimator for this type of experiments. This
estimator may be simply computed through standard algorithms for
logistic regression and does not require to formulate assumptions on
the distribution of the latent variables given the covariates. It
also has interesting asymptotic and finite-sample properties which
are maintained even with missing responses.

One of the basic assumptions on which the approach relies is that a
subject is assigned to the control arm or to the treatment arm with
a probability depending only on the observable covariates and not on
the pre-treatment response variable. Indeed, we could relax this
assumption, but we would have  more complex expressions for the
conditional probability of the response variables given their sum.
Similarly, the approach can be extended to the case in which
subjects assigned to both arm can access the treatment and then
non-compliance may also exist in the control arm, i.e. certain
subjects assigned to the control may decide instead to take the
treatment. The model presented in Section \ref{sec:ass} can be
easily extended to this case. Using the terminology of
\cite{Angrist:1996}, we have to consider the subpopulations of {\em
compliers}, {\em never-takers} and {\em always-takers}. By
exploiting an approximation similar to that illustrated in Section
\ref{sec:cond}, we can set up an estimator of the causal effect of
the treatment also in this case. The causal effect is again referred
to the subpopulation of compliers and is measured on the logit
scale. The approach would be complicated by the fact that the true
model involves a mixture on three subpopulations. Moreover, the
effect of the control on never-takers and that of the treatment on
compliers is not directly observable from the treatment arm as it
was in the original model. However, the resulting estimator would
maintain its simplicity as main advantage, being based on a series
of logistic regressions with suitable design matrices, which can be
performed by standard algorithms.

As a final comment consider that, driven by the application on the
BSE dataset, we only considered the case of repeated binary response
variables. However, the approach may be easily extended to the case
of response variables having a different nature (e.g. counting),
provided that the conditional distribution of these variables
belongs to the natural exponential family and the causal effect is
measured on a scale defined according to the canonical link function
for the adopted distribution \cite{McCullagh:1989}.
\section*{Acknowledgements}
I thank Dr. S. Ferro (Direzione Generale Sanit\`a e Politiche
Sociali, Regione Emilia Romagna, IT) for providing the dataset on
breast self examination and Dr. L. Grilli (University of Florence,
IT) for interesting discussions on the approach. I acknowledge the
financial support from Ministero dell'Istruzione,
dell'Universit\`{a} e della Ricerca (PRIN 2005 - ``Modelli marginali
per variabili categoriche con applicazioni all'analisi causale").

\section*{Appendix}
\subsection*{A1: Mathematical details on the approximation}
{\small In order to derive (\ref{eq:approx1}) consider that
\[
\log g(1,0|\b v)-\log g(0,1|\b v)\approx \log g_0(1,0|\b v)-\log
g_0(0,1|\b v)+b(\b v)\tr\b\be[h(0,1|\b v)-h(1,0|\b v)].
\]
Moreover, the first difference at right hand side is equal to $\b
a(\b v,0)\tr\b\al$ whereas
\[
h(0,1|\b v)=\frac{1}{g_0(1,0|\b v)e^{t(0,\bl
v,0)}}\int\frac{e^{\la(u,\bl v)}}{1+e^{\la(u,\bl
v)}}\frac{e^{t(0,\bl v,0)}}{[1+e^{\la(u,\bl v)+t(0,\bl
v,0)}]^2}\pi(1|u,\b v)\phi(u|\b v)du
\]
and then $h(0,1|\b v)-h(1,0|\b v)=h(\b v)$ as defined in
(\ref{eq:r}).
\subsection*{A2: Computation of $\hat{\b\Si}(\hat{\b\eta},\hat{\b\th})$}
{\small The derivative of $\b s(\b\eta,\b\th)$ has the following
structure
\[
\hat{\b
H}=\pmatrix{{\displaystyle\frac{\pa\ell_1(\hat{\b\eta})}{\pa\b\eta\pa\b\eta\tr}}
& \b O\cr
{\displaystyle\frac{\pa\ell_2(\hat{\b\eta},\hat{\b\th})}{\pa\b\th\pa\b\eta\tr}}
&
{\displaystyle\frac{\pa\ell_2(\hat{\b\eta},\hat{\b\th})}{\pa\b\th\pa\b\th\tr}}},
\]
where $\b O$ denotes a matrix of zeros of suitable dimension. The
first block of $\hat{\b H}$ corresponds to
\[
\frac{\pa\ell_1(\b\eta)}{\pa\b\eta\pa\b\eta\tr}=-\sum_i\frac{z_i}{p(z_i|\b
v_i)}\pi(0|\b v_i)\pi(1|\b v_i)\b m(\b v_i)\b m(\b v_i)\tr,
\]
whereas
\begin{eqnarray*}
\frac{\pa\ell_2(\b\eta,\b\th)}{\pa\b\th\pa\b\eta\tr}&=&\sum_id_i\bigg[y_{i2}-\frac{e^{\bl
w(\bl v_i,z_i,x_i)\tr\bl\th}}{1+e^{\bl w(\bl
v_i,z_i,x_i)\tr\bl\th}}\bigg]\frac{\pa\b w(\bl
v_i,z_i,x_i)}{\pa\b\eta\tr}+\\&-&\sum_id_i\frac{e^{\bl w(\bl
v_i,z_i,x_i)\tr\bl\th}}{[1+e^{\bl w(\bl v_i,z_i,x_i)\tr\bl\th}]^2}\b
w(\bl v_i,z_i,x_i)\b w(\bl v_i,z_i,x_i)\tr\frac{\pa\b w(\bl
v_i,z_i,x_i)}{\pa\b\eta\tr},
\end{eqnarray*}
with
\[
\frac{\pa\b w(\bl v_i,z_i,x_i)}{\pa\b\eta\tr}=\pmatrix{\b 0\cr
(1-z_i)\pi(0|\b v_i)\pi(1|\b v_i)\b b(\b v_i)}\b m(\b v_i)\tr.
\]
Finally, we have
\[
\frac{\pa\ell_2(\b\eta,\b\th)}{\pa\b\th\pa\b\th\tr}=-\sum_id_i
\frac{e^{\bl w(\bl v_i,z_i,x_i)\tr\bl\th}}{[1+e^{\bl w(\bl
v_i,z_i,x_i)\tr\bl\th}]^2}\b w(\bl v_i,z_i,x_i)\b w(\bl
v_i,z_i,x_i)\tr.
\]

The estimate of the variance-covariance matrix of the score may be
expressed as
\[
\hat{\b K} = \pmatrix{\hat{\b K}_{11}& \hat{\b K}_{12}\cr \hat{\b
K}_{21}& \hat{\b K}_{22}},
\]
where
\begin{eqnarray*}
\hat{\b K}_{11} &=& \sum_i
\frac{z_i}{p(z_i|\b v_i)^2}[x_i-\hat{\pi}(1|\b v_i)]^2\b m(\b v_i)\b m(\b v_i)\tr,\\
\hat{\b K}_{12} &=& \sum_i \frac{d_iz_i}{p(z_i|\b
v_i)}[x_i-\hat{\pi}(1|\b v_i)]\bigg[y_{i2}-\frac{e^{\bl w(\bl
v_i,z_i,x_i)\tr\bl\th}}{1+e^{\bl w(\bl
v_i,z_i,x_i)\tr\bl\th}}\bigg]\b m(\b v_i)\b w(\b v_i,z_i,x_i)\tr,\\
\hat{\b K}_{22} &=& \sum_i d_i\bigg[y_{i2}-\frac{e^{\bl w(\bl
v_i,z_i,x_i)\tr\bl\th}}{1+e^{\bl w(\bl
v_i,z_i,x_i)\tr\bl\th}}\bigg]^2\b w(\bl v_i,z_i,x_i)\b w(\bl
v_i,z_i,x_i)\tr
\end{eqnarray*}
and $\hat{\b K}_{21}=\hat{\b K}_{12}\tr$.}
\end{document}